# L'utilità di una teoria inutile
# Crittografia, firma digitale e teoria dei numeri[1]

**Alfredo Esposito[2]**
**InfoCert S.p.A.**

**February 6, 2011**

**Riassunto** La teoria dei numeri era considerata la parte della matematica meno utile alla vita quotidiana. La crittografia invece era vista come qualcosa riservata ai militari e ai diplomatici. Qui si mostra come le due cose siano oggi strettamente collegate e come questo legame possa avere influenza sulla nostra vita quotidiana.

**Abstract** The theory of numbers was supposed to be the less useful branch of mathematics. At the same time, cryptography was thought to be a military or a diplomatic issue. In this note we show how the two concepts are today strictly related and how this link affects our daily activities.

## 1. Introduzione

Nel suo famoso libro *Apologia di un matematico* [5], G.H Hardy, noto per i suoi contributi alla teoria dei numeri ed all'analisi e per avere intuito le incredibili potenzialità del giovane talento indiano Srinivasa Ramanujan, scriveva[3] "*Nessuno ha ancora scoperto un uso bellico della teoria dei numeri o della relatività...*" e si vantava della "inutilità" pratica della matematica di cui si occupava, il cui unico valore era estetico.

Quando morì, nel 1947, l'applicazione bellica della relatività[4] era stata purtroppo vista all'opera su Hiroshima e Nagasaki. Fosse vissuto più a lungo, parecchio più a lungo, avrebbe scoperto, forse con sgomento, che anche la teoria dei numeri, il problema della distribuzione dei numeri primi e il problema della fattorizzazione di un numero composto grande, hanno un utilizzo pratico, bellico e non solo. La troviamo, infatti, negli acquisti che facciamo su Internet, nelle applicazioni di home banking, nella posta elettronica certificata, nella possibilità di lavorare in remoto anche su dati e sistemi che richiedono riservatezza, fino alla firma digitale, l'equivalente legale, nel mondo digitale, della firma autografa.

---

[2] alfredo.esposito@infocert.it, alfespo@yahoo.it

[3] pag. 99.

[4] Segnatamente l'equivalenza massa-energia $E=mc^2$

Senza avere la pretesa di affrontare i dettagli più minuti, l'articolo vuole fornire gli elementi di base per comprendere come l'*inutile* teoria dei numeri sia diventata una componente importante della vita quotidiana di chiunque.

Inizieremo parlando brevemente della crittografia e dei suoi scopi, dandone un ridottissimo *excursus* storico e qualche accenno della sua evoluzione da arte a scienza. A seguire descriveremo i due algoritmi più usati nella crittografia asimmetrica, che è la protagonista principale di questa nota. Nella parte finale descriveremo l'applicazione pratica dei concetti illustrati, con enfasi particolare sulla firma digitale, che ha permesso di sostituire completamente la carta in parecchi ambiti. Infine torneremo brevemente sui concetti generali, discutendo la ‘robustezza’ degli algoritmi crittografici.

Tre appendici completano il testo. La prima fornisce la matematica di base, che ho preferito risparmiare al lettore già in possesso del (minimo) bagaglio matematico necessario,[5] la seconda illustra i concetti base relativi alla sicurezza delle informazioni, utili a comprendere meglio alcuni passaggi del testo, ma non indispensabili ad una prima lettura, la terza raggruppa le definizioni di alcuni termini tecnici.

La bibliografia, infine, integra le referenze date nel testo con ulteriori suggerimenti di lettura.

## 2. La crittografia e la steganografia

I requisiti di *integrità*, *riservatezza* e *non ripudio*[6] non sono certo un’invenzione dell’era dell’informatica.

Fin dall’antichità, essenzialmente per motivi militari, si è posto il problema di garantire l’*integrità* dei messaggi, utilizzando sigilli di varia natura, il *non ripudio* ovvero l’autenticità del mittente (con la firma autografa oppure apponendo simboli sui sigilli stessi, per esempio imprimendo il simbolo con un anello) e la *riservatezza*, tramite la crittografia e la steganografia.

Per essere più precisi, la steganografia è l'arte di nascondere un'informazione in un testo che, per altri versi, appare perfettamente normale; un esempio banalissimo è dato dalla filastrocca che aiutava, quando l'autore di queste note era assai più giovane, a ricordare, tramite le iniziali delle parole, la sequenza ovest-est delle Alpi.[7] In questo caso non c'era riservatezza, non c'era nulla da nascondere, ma il metodo è identico: il testo della filastrocca racchiude un altro testo (le iniziali). In maniera più seria, un'informazione può essere nascosta in un'immagine alterando i bit meno significativi di ciascun pixel, in maniera che neanche un occhio allenato possa cogliere l'alterazione.

---

[5] Rimanendo comunque all’algebra liceale, evitando, a spese del rigore, la terminologia più precisa dell’algebra astratta.

[6] Cfr. Appendice 2.

[7] “MACOnGRAnPENaLERECAGIU” dove le prime sillabe di ogni settore alpino (MArittime, COzie, GRAie, PENnine, LEpontine, REtiche, CArniche e GIUlie) concorrono a formare un testo mnemonico.

La crittografia non si cura invece di dare al messaggio un aspetto innocuo. Chi intercetta un messaggio crittografato sa che là dentro ci sono informazioni che potrebbero essere utili. Tuttavia, se non è in possesso della chiave giusta, non avrà modo di accedervi.

## 2.1. Un po' di storia

Della steganografia non parleremo oltre, salvo citare come esempi gli inchiostri cosiddetti simpatici (che svaniscono e riappaiono se riscaldati) e la microstampa (*microdot*) usata da agenti tedeschi durante la II Guerra Mondiale e spesso presente nei fumetti, con le clausole più vessatorie scritte nei puntini delle ‘i’ dei contratti capestro che Paperino doveva firmare con Zio Paperone.

Uno degli esempi più antichi di crittografia di cui si è a conoscenza, è lo scitale degli spartani (*scitala lacedemonica*). Questo era un bastone di legno attorno al quale era avvolta una striscia di stoffa e/o pergamena. Il messaggio era scritto mentre la pergamena era avvolta sul bastone; quindi la pergamena veniva svolta e la sequenza dei caratteri appariva casuale. Solo riavvolgendola su un bastone dello stesso spessore, il messaggio sarebbe stato nuovamente leggibile. La pergamena poteva inoltre essere utilizzata come cinta (una forma di steganografia).

Questi sistemi si definiscono ‘crittografia per trasposizione’ perché i simboli che compongono il messaggio crittografato (*ciphertext*), oltre a provenire dallo stesso alfabeto, sono esattamente quelli presenti[8] nel testo originario (*plaintext*).

In alternativa, la crittografia può essere ‘per sostituzione’: ogni simbolo o gruppo di simboli è sostituito da un altro, che può o meno far parte dello stesso alfabeto. Un celebre esempio storico è quello descritto da Giulio Cesare nel *De bello gallico*:[9] un messaggio in cui i caratteri latini erano stati sostituiti da caratteri greci, sfruttando l'ignoranza, reale o presunta, dei Galli.

Lo stesso Cesare, almeno a quanto ci racconta Svetonio,[10] usava un altro algoritmo di sostituzione, sostituendo ogni lettera del messaggio con quella spostata di tre passi nel medesimo alfabeto latino. Per esempio, se il *plaintext* è

*Nel mezzo del cammin di nostra vita*

una trasposizione di tre caratteri a destra, assumendo l’alfabeto corrente a 26 caratteri e che dopo la Z torni la A, produce il *ciphertext*

*Qho phccr gho fdpplq gl qrvwud ylwd.*

Un altro esempio storico, rilevante per le conseguenze[11] che ebbe sulla storia

---

[8] Il testo cifrato è quindi una permutazione del testo originario.

[9] Libro V, 48.

[10] Svetonio, *De vita caesarum*.

[11] La decifrazione condusse Maria Stuarda al patibolo e durante il regno di Elisabetta furono poste le basi della futura potenza commerciale e marittima della nazione ed iniziò la colo-

dell'Inghilterra e, di riflesso, dell'Europa, fu la decifrazione, narrata nel libro di Singh *Codici e segreti* [14],[12] del codice usato da Maria Stuarda con i suoi sodali nella lotta contro Elisabetta I.

Per molto tempo la segretezza dei messaggi è stata legata esclusivamente alla segretezza dell'algoritmo di cifratura o alla condivisione di apparati meccanici o, più recentemente, elettromeccanici installati presso tutti i possibili partner. Il più noto, il tedesco Enigma, è stato protagonisti di film e libri[13] e la sua decifratura, cui diede un contributo fondamentale Alan Turing [7], fu un fattore decisivo per vincere la guerra nell'Oceano Atlantico e, forse, per la sconfitta del nazifascismo.

Nel 1800, pur rimanendo ancora un po' artigianale, inizia un approccio più matematico alla crittografia.

Uno dei risultati più rilevanti è il principio di Kerckhoffs, riformulato indipendentemente da Shannon, che afferma che la forza di un sistema crittografico non deve risiedere nella segretezza dell'algoritmo, ma della chiave. Questo approccio è tuttora oggetto di dibattito, ma la tesi che rendere pubblico un algoritmo renda più probabile la scoperta delle sue, eventuali, debolezze sembra piuttosto convincente.[14]

All'inizio del Novecento viene messo a punto l'unico sistema crittografico di cui è possibile (finora) dimostrare matematicamente l'inviolabilità (Vernam, One-Time-Pad); risulta però piuttosto complicato da adottare, richiedendo una nuova chiave per ogni messaggio, ciascuna lunga quanto il messaggio stesso.

Dopo la seconda guerra mondiale, i lavori di Shannon [13] stabilirono una base teorica per la crittografia e la crittoanalisi, trasformandole da arte in scienza.

Fino alla metà degli anni '70 la crittografia rimase quasi di esclusiva competenza dei circoli militari; a metà del decennio, lo statunitense National Bureau of Standards (oggi NIST, National Institute for Standard and Technology [27]) pubblicò il primo standard crittografico disponibile per le applicazioni commerciali: il DES, che dal 2001 è stato sostituito da un nuovo standard (AES) selezionato dopo un concorso pubblico che ha coinvolto i crittografi di tutto il mondo.

## 2.2. La crittografia simmetrica (a chiave segreta)

In tutti i metodi usati fino al 1976, per garantire la riservatezza dei messaggi gli interlocutori dovevano preventivamente condividere un segreto, fosse esso un metodo o un algoritmo oppure la chiave di cifratura utilizzata dall'algoritmo stesso.

La figura 1 illustra come opera la crittografia a chiave segreta ed evidenzia

---

nizzazione dell'America settentrionale.

[12] Capitolo 1

[13] http://ed-thelen.org/comp-hist/NSA-Enigma.html, visitato il 10/03/2010

[14] Non diversamente da quanto accade normalmente nell'attività scientifica, dove le pubblicazioni sono riviste da esperti e gli esperimenti devono essere descritti e ripetibili.

perché sia definita simmetrica: la chiave di cifratura e decifratura sono identiche (K) ed il flusso può andare da destra a sinistra o viceversa.

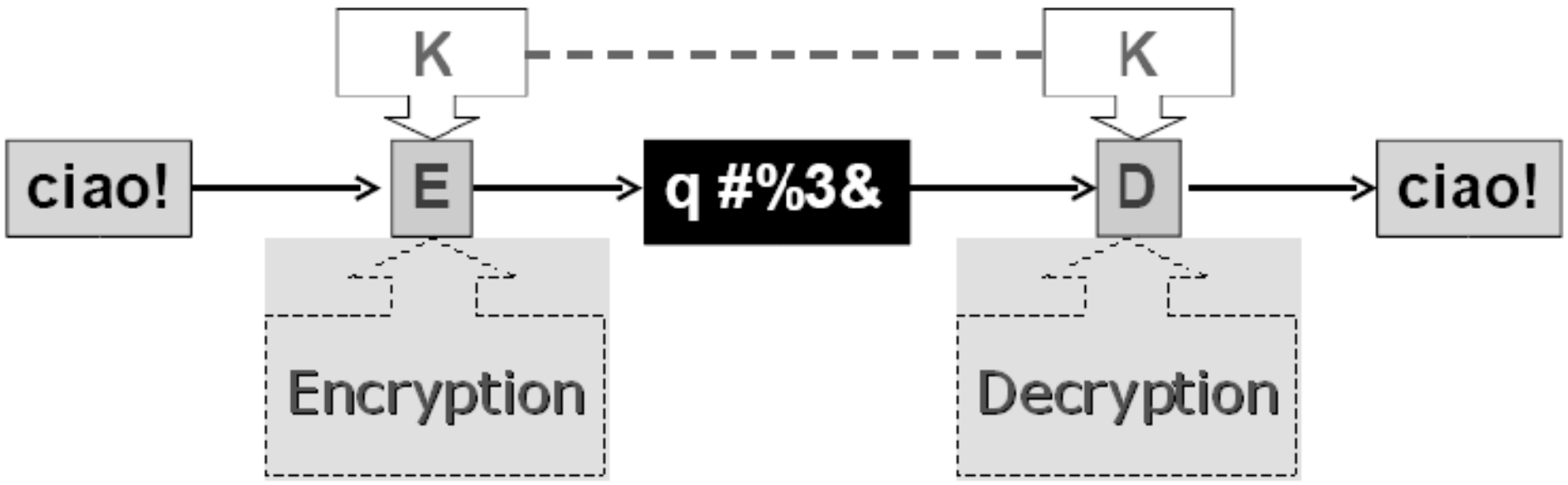

Figura 1 – La crittografia simmetrica

L'uso di tale tipo di crittografia era reso possibile (e lo è ancora in certi contesti) dal fatto che la comunità degli interlocutori era limitata e nota a priori, permettendo la diffusione della chiave (il *segreto condiviso*) secondo meccanismi opportuni (per esempio incontrandosi fisicamente) distinti da quelli utilizzati poi nello scambio di messaggi.

Se gli interlocutori diventano numerosi il problema della distribuzione delle chiavi diventa oneroso e difficile da gestire.

Immaginiamo infatti di avere 10 persone che hanno deciso di tenere riservate le loro comunicazioni.

Per comodità chiamiamole $P_i$ dove il suffisso *i* può valere 1,2,...10. Tutte le conversazioni sono crittografate in modo che nessuno degli altri possa carpirne il contenuto.

Allora $P_1$ deve incontrarsi con $P_2$, $P_3$, $P_4$, $P_5$, $P_6$, $P_7$, $P_8$, $P_9$, $P_{10}$ e gestire quindi 9 chiavi diverse.

A sua volta $P_2$ deve incontrarsi con $P_3$ ...fino a $P_{10}$ e qui sono 8 chiavi.

Senza ripetere la storia per tutti i $P_i$, abbiamo che le chiavi in ballo sono

**9+8+7+6+5+4+3+2+1=45**

Già su un campione così limitato, le chiavi in ballo sono ben 45.

In generale se le persone sono *N*, le chiavi da gestire sono le combinazioni di *N* oggetti presi due per volta

$$N\frac{(N-1)}{2} \cong N^2 \quad per\ N >> 1 \tag{1}$$

E' palese che si tratta di una situazione che diventa rapidamente ingestibile al crescere di *N*. Ovviamente, in alcuni casi reali le cose sono più semplici; per esempio, se gli interlocutori fanno tutti riferimento ad un'unica Amministrazione Centrale che funge da snodo e parlano solo con questa e non tra loro, allora lo snodo centrale distribuisce le chiavi agli uffici periferici ed il numero di chiavi da gestire cresce come *N*. Resta comunque da risolvere il problema della distribuzione *a priori* della chiave su un canale sicuro, prima di poter avviare lo scambio ‘sicuro’ elettronico.

Ma come ottenere lo stesso risultato qualora gli interlocutori non abbiano alcun modo di incontrarsi?

### 2.3. La crittografia asimmetrica (a chiavi pubbliche)

Nel 1976 R. Merkle[15] e, indipendentemente, W. Diffie e M. Hellman [3] inventarono il sistema di crittografia a chiavi pubbliche (Public Key Cryptography - PKC), che risolve in un sol colpo il problema della diffusione del segreto condiviso, dell'integrità del messaggio e del non ripudio.

Si è trattato di una vera rivoluzione. Duemila e più anni di storia della crittografia, fosse essa gestita con carta e penna o con i più moderni sistemi automatici, venivano rivoluzionati da questa idea che risolveva il problema della distribuzione delle chiavi di cifratura.[16]

La paternità della scoperta/invenzione è stata successivamente rivendicata dai servizi di *intelligence* britannici. Il CESG – Communications Electronics Security Group (allora si chiamava GCHQ – Government Communications Headquarters) afferma che, prima degli autori citati, la tecnica sia stata inventata da Ellis, un matematico che operava presso di loro, nei primi anni '70, ma tenuta segreta per motivi di sicurezza nazionale.[17]
Vera o falsa che sia questa rivendicazione, la crittografia asimmetrica comincia la sua 'carriera' effettiva con la pubblicazione del lavoro di Diffie-Hellman, basato sulla difficoltà di invertire il calcolo dei logaritmi discreti[18], non possedendo tutte le informazioni su come sono stati calcolati.

Lo schema si basa sull'utilizzo di due chiavi, una per cifrare ed una per decifrare. La proprietà fondamentale di tutti i sistemi di crittografia a chiavi pubbliche è che le due chiavi sono legate da una relazione matematica per cui un testo manipolato con una delle due può essere decifrato solo con l'altra, ma le chiavi sono diverse, da cui l'aggettivo asimmetrica.

Nella PKC, una delle chiavi è resa *pubblica* (da cui il nome) mentre l'altra (detta *privata*) non deve essere mai rivelata.

Un mittente, Alice,[19] che voglia inviare un messaggio riservato al destinatario Bob, lo cifra con la chiave pubblica del destinatario Bob, il quale, essendo il possessore della chiave privata corrispondente, sarà l'unico in grado di leggerlo.[20]

Analogamente se il mittente Alice cifra il messaggio con la propria chiave privata, chiunque potrà aprire il messaggio e verificarne il mittente, utilizzando la chiave pubblica di Alice (*non ripudio*).

---

15 http://www.merkle.com/1974/ visitato il 09/06/2010.

16 La possibilità di applicare questo tipo di crittografia non è ovviamente indipendente dalla disponibilità di computer in grado di svolgere in tempi brevi calcoli aritmetici complessi.

17 Sono in linea alcune delle note originali (http://www.cesg.gov.uk/publications/media/notense.pdf, visitato il 21/05/2010).

18 A valori interi.

19 Convenzionalmente i due interlocutori di uno scambio crittografato sono chiamati Alice e Bob, a cui si aggiunge spesso Eve (da eavesdropper, che cerca di intercettare i dati sulla rete).

20 Neanche Alice sarebbe in grado di farlo!

Per fare un paragone con il mondo degli oggetti più consueti, immaginate di consegnare a mano un messaggio, imbucandolo direttamente nella cassetta delle lettere del destinatario e supponiamo che la cassetta sia blindata e che, ovviamente, la chiave (o la combinazione) della cassetta sia saldamente nelle mani di Bob. Neanche Alice può riprenderlo una volta che lo ha inviato. Oppure immaginate che Bob metta a disposizione in una bacheca centinaia di lucchetti tutti uguali. Alice prende il messaggio (o l'oggetto) che vuole inviare a Bob, lo inserisce in una scatola di metallo e la chiude con uno dei lucchetti. Solo Bob ha la chiave per aprire la scatola ed accedere al contenuto.

Virtualizzate il tutto: il lucchetto è uno solo (la chiave pubblica), ma copiabile quante volte volete, ed avete la cifratura tramite la crittografia asimmetrica. Trovare un'analogia per il *non ripudio* è un po' più complicato: la situazione che potete provare ad immaginare è quella, peraltro già implementata da alcune banche, di una copia elettronica dello *specimen*, cioè dell'immagine della firma autografa, depositata in ogni agenzia. Ogni agenzia è quindi in grado di verificare la vostra firma, indipendentemente dal fatto che voi ci siate mai stati fisicamente. Questo paragone evidenzia solo un aspetto della firma elettronica (non dovete incontrarvi fisicamente e preliminarmente con l'interlocutore); per tutti gli altri versi, la firma elettronica è più sicura e lega in maniera assai più forte la firma al firmatario, in quanto, come vedremo, anche la minima modifica del contenuto firmato (anche un singolo bit!) farebbe fallire l'operazione di verifica.

Il principio che sta quindi alla base di questa crittografia è l'esistenza di algoritmi o procedure che abbiano tre caratteristiche fondamentali:

La procedura deve essere:

- facile[21] da eseguire affinché tutti possano scrivere a Bob
- facile da invertire, cioè da verificare/decifrare, <u>possedendo</u> l'informazione giusta (la chiave privata)
- molto difficile da invertire, cioè da verificare/decifrare <u>non possedendo</u> l'informazione giusta.

## 3. Crittografia e complessità computazionale

L'esistenza della crittografia a chiave pubblica è legata strettamente all'esistenza di funzioni monodirezionali (*one-way function*) cioè di funzioni che siano facilmente calcolabili in un verso, ma per le quali sia praticamente impossibile risalire dal risultato al valore originario.

In maniera più formale, dati due insiemi $X$ e $Y$ , la funzione $f$ è monodirezionale se

[21] Nel seguito proveremo a dare agli aggettivi ‘facile’ e ‘difficile’ un significato più preciso nel nostro contesto. Per ora è comunque sufficiente il significato che hanno nel linguaggio comune.

$$f(x): X \rightarrow Y \text{ è tale che } \forall\, x \in X \text{ c'è un corrispondente } y \in Y = f(x) \text{ e}$$
$$\forall y \in Y \text{ è } \textit{praticamente} \text{ impossibile[22] trovare la } x \text{ tale che } f(x) = y$$

Per darvi un idea di una tipica operazione *one-way*, pensate agli elenchi telefonici tradizionali. Avete un cognome e dovete cercare il numero. Soluzione pressoché immediata, anche su un grosso volume. Ma se avete solo il numero e dovete cercare l'abbonato corrispondente? Togliendo i numeri speciali riportati nelle prime pagine, dovreste scorrere tutto l'elenco. Il numero atteso di numeri da scorrere è $N/2$, dove $N$ è il numero di voci che costituiscono l'intero elenco.

Una funzione *one-way with trapdoor* è una funzione che ha la caratteristica ulteriore che, in presenza di un'informazione aggiuntiva, l'inversione torna ad essere facile.

L'esistenza o meno delle funzioni *one-way* si ricollega al tema della complessità computazionale,[23] ovvero a quella branca della matematica (o dell'informatica teorica) che studia la difficoltà di risoluzione di un problema tramite un algoritmo implementabile su un computer.[24] In generale, è possibile assegnare ad ogni problema un numero $n$ che è rappresentativo della 'dimensione' del problema: se stiamo ordinando degli oggetti sarà il loro numero, se stiamo cercando l'abbonato dell'esempio sarà il numero di abbonati presenti sull'elenco, se stiamo parlando di funzioni crittografiche sarà la dimensione in bit della chiave adoperata. La modalità con cui il tempo[25] di risoluzione del problema (in termini non di secondi, ma di passi elementari di programma) viene a dipendere da $n$ caratterizza la complessità del problema.

La complessità algoritmica è un argomento dove si incrociano matematica, filosofia, informatica teorica ed è un campo dove è tuttora aperto uno dei Millennium problems, per i quali il Clay Mathematics Institute[26] promette un premio di 1 milione di dollari per il primo che riuscirà a risolverli.

Il problema in questione è espresso nei seguenti termini:

$$\text{E' vera la relazione } P \neq NP\,?$$

Vediamo di spiegarla. I simboli $P$ e $NP$ indicano due differenti raggruppamenti di problemi, sulla base della loro 'facilità' di risoluzione.

Con classe di complessità P si intende l'insieme dei problemi che possono essere risolti in tempo polinomiale, cioè tali che il tempo di risoluzione è propor-zionale a $n^a$ dove $a$ è una qualche costante. Un problema che ammetta una soluzione di questo tipo è 'facile' da risolvere, perché il lavoro da svolgere cresce lentamente (relativamente) con la complessità del problema.

---

[22] E' ovviamente possibile se $Y$ è un insieme finito, ma il costo elaborativo è elevatissimo e lunghissimo il tempo di elaborazione necessario.

[23] Il primo capitolo di Goldreich *Computational Complexity A Conceptual Perspective* (disponibile in rete http://www.wisdom.weizmann.ac.il/~oded/books.html) è un'ottima introduzione all'argomento (visitato il 21/06/2010).

[24] Ci si riferisce ad una macchina di Turing, di cui i computer attuali sono una realizzazione pratica.

[25] Il concetto è in realtà più ampio e sulla complessità interviene anche la quantità di memoria necessaria al programma per memorizzare i dati intermedi.

[26] http://www.claymath.org/millennium/P_vs_NP/ visitato il 09/06/2010.

Esistono altri problemi (esponenziali) per i quali il tempo di risoluzione cresce invece come $a^n$ quindi in maniera più esplosiva.

La tabella seguente esemplifica la differenza per $a=2$.

*Tabella 1*

| $n$ | $n^2$ | $2^n$ |
|---:|---:|---:|
| 1 | 1 | 2 |
| 5 | 25 | 32 |
| 10 | 100 | 1024 |
| 15 | 225 | 32768 |
| 20 | 400 | 1048576 |

I problemi classificati come *NP* sono invece i problemi che hanno una procedura rapida di verifica delle soluzioni, quando se ne ha una a disposizione, ma per i quali non è detto che sia facile trovare tale soluzione (*N* sta per *non deterministico*).

Ovviamente $P \subset NP$, ma la domanda che i matematici si fanno è se le due classi sono effettivamente diverse l'una dall'altra. Se si riuscisse a dimostrare che $P \neq NP$ allora potremmo essere certi che una funzione *one-way* esiste davvero. In mancanza di questa dimostrazione, dobbiamo basarci sulle congetture e usare per la nostra crittografia asimmetrica problemi che sono 'generalmente ritenuti' *NP*. Problemi di questo tipo e che incontreremo tra poco sono la fattorizzazione di un numero composto[27] e il calcolo del logaritmo discreto.

## 4. Gli algoritmi della crittografia a chiave pubblica

E' il momento di entrare nel dettaglio degli algoritmi che realizzano la crittografia a chiave pubblica. Ne vedremo due, di gran lunga i più usati, che si appoggiano a due diversi problemi 'difficili':

- RSA, basato sulla difficoltà di fattorizzare un numero grande;
- Diffie-Hellman, che si basa invece sulla difficoltà di calcolare il logaritmo discreto.

Nel seguito faremo uso delle congruenze e dell'aritmetica modulare, di cui

[27] Cfr. Appendice 1.

l’appendice 1 fornisce i concetti fondamentali.

### 4.1. L'algoritmo RSA

Inventato nel 1978 da Rivest, Shamir e Adleman [12], dalle cui iniziali prende il nome, è il sistema di crittografia a chiavi pubbliche attualmente più diffuso ed utilizzato.

L'RSA si basa essenzialmente sulla difficoltà di scomporre un numero grande nei suoi fattori primi. L'equivalenza formale della fattorizzazione e dell'RSA non è stata dimostrata, quindi non è da escludere (per quanto improbabile) che qualcuno possa trovare un modo di violare l'RSA anche se il problema della fattorizzazione rimanesse irrisolto. Tuttavia lo stretto legame esistente lo si può dedurre dal fatto che, al contrario, la soluzione al problema della fattorizzazione implica automaticamente la soluzione dell'RSA.

Per iniziare occorre scegliere due interi grandi $p$ e $q$ che siano primi. Se ne calcola il prodotto $N = p \cdot q$ che sarà il *modulo* adottato per tutte le operazioni susseguenti; $p$ e $q$ devono essere tenuti segreti. Si sceglie quindi un altro intero positivo $e$, che sia relativamente primo rispetto al prodotto $\varphi(N) = (p-1)\cdot(q-1)$ cioè tale che il massimo comune divisore tra $e$ e $\varphi$ sia 1, con $1 < e < \varphi$. La coppia $\{N, e\}$ rappresenta la chiave pubblica.

Si calcola quindi un intero $d$ tale che

$$e \cdot d \equiv 1 \ mod \ \varphi(N) \qquad (2)$$

La coppia $\{N, d\}$ è la chiave privata.

A questo punto per cifrare un messaggio $M$, lo si divide in blocchi $M_i$ di lunghezza pari a $N$ (completando opportunamente l’ultimo blocco) e si calcola il messaggio cifrato (*ciphertext*).

$$C_i \equiv (M_i)^e \ mod \ N \qquad (3)$$

Per riottenere il messaggio originale si divide il *ciphertext* in blocchi di lunghezza $N$ e si calcola

$$M_i \equiv (C_i)^d \ mod \ N \qquad (4)$$

Essendo infatti

$$(C_i)^d \equiv (M_i^{\,e})^d \ mod \ N \qquad (5)$$

si ha, sfruttando la (2) e le proprietà di riduzione degli esponenti proprie delle congruenze,

$$(M_i^{\,e})^d \equiv (M_i)^{e \cdot d \ mod \ \varphi(N)} = M_i \qquad (6)$$

La conoscenza di $p$ e $q$, e quindi di ($p$-1) e ($q$-1), permetterebbe di ottenere $d$ piuttosto facilmente, risolvendo la (2) ed essendo noto l’esponente $e$: quindi la sicurezza dell’algoritmo RSA dipende dalla difficoltà, nel senso definito prima, di determinare i fattori primi di un numero intero grande. Inversamente, si può dimostrare che la conoscenza di $d$ o di $\varphi(N)$ fornisce un metodo efficiente per la fattorizzazione di $N$.

### 4.1.1. Un esempio di cifratura RSA

Per fare un esempio banale, scegliamo $p$=19 e $q$=17. Il modulo risultante è $N$=323, mentre $\varphi$=288. Per semplicità di calcolo, il testo che cifreremo sarà composto da un singolo carattere, la lettera C. Scegliamo ancora un semplice $e$=17, da cui deriva $d$=17, il che rende le potenze relativamente maneggevoli.

Consideriamo una codifica ad hoc dei caratteri, ammettendo solo le 26 lettere dell'alfabeto, codificate numericamente con la loro posizione.

Quindi A→1, B →2, C→3, e così via fino a Z→26.

La codifica numerica della nostra C è quindi 3.

Elevando all'esponente pubblico 17, operazione alla portata di una calcolatrice, si ottiene il valore decimale 241, calcolato modulo 323, che rappresenta il nostro messaggio cifrato.

Facciamo ora l'operazione inversa. Per farlo dobbiamo sfruttare le proprietà delle congruenze per poter maneggiare le potenze.[28]

Per decifrare il messaggio dobbiamo calcolare:

$$241^{17}\ (mod\ 323)$$

Possiamo scrivere $241^{17}$ come $241^3 \cdot 241^3 \cdot 241^3 \cdot 241^3 \cdot 241^3 \cdot 241^2$ e considerando che $241^3 \equiv 316\ (mod\ 323)$ e che $241^2 \equiv 264\ (mod\ 323)$, il nostro calcolo diventa

$$316 \cdot 316 \cdot 316 \cdot 316 \cdot 316 \cdot 264\ (mod\ 323)$$

cioè

$$(-7) \cdot (-7) \cdot (-7) \cdot (-7) \cdot (-7) \cdot 264\ (mod\ 323)$$

Notando che

$$7^3 = 343 \equiv 20\ (mod\ 323)$$

e che

$$(-7 \cdot 264) = -1848 \equiv 90\ (mod\ 323)$$

abbiamo

$$(-20) \cdot (-7) \cdot 90\ \equiv 3\ (mod\ 323)$$

e infine ritroviamo la nostra C.

I parametri di questo esempio sono stati scelti opportunamente per rendere fattibili i calcoli semi manuali. Nella realtà si usano numeri assai più grandi. Le chiavi RSA usate comunemente sono oggi lunghe almeno 1024 bit, corrispondenti a numeri dell'ordine di $2^{1024}$.

## 4.2. L'algoritmo di Diffie-Hellman

L'altro algoritmo largamente usato di crittografia asimmetrica è quello di Diffie-Hellman, cioè quello che nel 1976 ha dato il via alla crittografia a chiave pubblica.

[28] cfr. Appendice 1.

Per accordarsi su una chiave segreta comune (simmetrica), Alice e Bob devono scegliere un numero primo $p$ (grande) e un numero naturale $g$ qualsiasi, $2 < g < p-2$ . Non è necessario che $p$ e $g$ siano tenuti segreti.

Ciascun partner seleziona poi un numero segreto (rispettivamente ***a*** e ***b***); i due calcolano i due valori

Alice calcola $$\alpha \equiv g^{a} \ (mod\ p) \tag{7}$$

Bob calcola $$\beta \equiv g^{b} \ (mod\ p) \tag{8}$$

I due partner si scambiano i valori così ottenuti.

A questo punto ciascuno dei due è in grado di calcolare il segreto comune $S$. Vediamo come, mostrando in parallelo i calcoli dei nostri due protagonisti.

| I calcoli di Alice (che conosce $a$) | …e quelli di Bob (che conosce $b$) |
|---|---|
| $\beta^{a} \equiv \left(g^{b}\right)^{a} \ (mod\ p)$ | $\alpha^{b} \equiv \left(g^{a}\right)^{b} \ (mod\ p)$ |

In questo modo Alice e Bob condividono una chiave segreta di cifratura

$$S = g^{(a \cdot b)} \tag{9}$$

con cui si garantiranno la riservatezza dei messaggi successivi. Eve, d'altra parte, dovrebbe calcolare $a$ oppure $b$ partendo dai messaggi che sono transitati in chiaro sulla rete ($\alpha$ o $\beta$). Per far questo dovrebbe calcolare il logaritmo in una base intera ($g$) conoscendo esclusivamente il resto della divisione per $p$, un problema ritenuto difficile nel senso descritto in precedenza.

Questa matematica, dall'aspetto vagamente esoterico, la usate tutte le volte che vi collegate al sito della vostra banca, perché è alla base del funzionamento dell'SSL,[29] cioè del protocollo che vi permette di identificare con ragionevole certezza che state parlando proprio con il sito della banca e di trasmettere sulla rete dati come la password o disposizioni bancarie senza temere che esse vengano intercettate.

## 5. Dalla teoria alla pratica

La bella matematica vista finora, quando deve essere tradotta in pratica, deve superare alcune difficoltà operative.

Il principale ostacolo è il tempo di calcolo necessario per eseguire certe operazioni matematiche. Le operazioni di elevamento a potenza nella crittografia asimmetrica sono di qualche ordine di grandezza più ‘pesanti’ di quelle relative alla crittografia simmetrica (a parità di livello di sicurezza).

Di conseguenza nel trasformare gli algoritmi da formule sulla carta a codice eseguibile, si è adottata una strategia mista, in cui la crittografia asimmetrica è utilizzata su sequenze relativamente corte, dell’ordine di qualche centinaio di bit, mentre il lavoro più pesante è svolto dai più rapidi algoritmi simmetrici, nel caso della cifratura, o dalle funzioni di *hash*, che vedremo tra poco, per il non ripudio.

[29] Il lucchetto chiuso che il vostro browser vi mostra quando vi collegate ad un sito https, per es. il sito della banca, indica che l’SSL è attivo.

### 5.1. La cifratura

L'invio di un messaggio cifrato segue quindi un iter leggermente più articolato. Definiamo due funzioni, *Enc (M, K)* che indica l'operazione di cifratura del messaggio *M* con la chiave *K* e *Dec* ($M_c$, *K)* che è la funzione inversa di decifratura del ciphertext $M_c$.

Il mittente Alice:

1. genera una chiave di cifratura <u>simmetrica</u> in maniera casuale (detta chiave di sessione *K*);
2. cifra il messaggio con la chiave di cifratura così generata: $M_c=Enc\ (M, K)$;
3. cifra la chiave di sessione[30] con la chiave pubblica del destinatario Bob ($B_{pk}$): $K_c=Enc\ (K, B_{pk})$;
4. invia il messaggio cifrato $M_c$ e la chiave di sessione cifrata $K_c$.

Il destinatario Bob:

1. riceve il messaggio composto da $M_c$ e $K_c$ ;
2. decifra con la propria chiave privata $B_{sk}$ la chiave di cifratura (chiave di sessione) $K=Dec\ (K_c,B_{sk})$ ;
3. decifra quindi il messaggio $M=Dec\ (M_c, K)$.

### 5.2. La 'digital signature'

Innanzitutto vorrei chiarire che questo paragrafo ha un titolo in inglese non per un vezzo di chi scrive, ma per evitare incomprensioni e confusioni tra la firma digitale dal punto di vista strettamente tecnico e la 'firma digitale' definita nell'ordinamento italiano [2], e che pure si basa esattamente sulla stessa tecnologia.

Dal punto di vista tecnico/crittografico, *digital signature* è l'applicazione della coppia di chiavi della crittografia asimmetrica invertendone il ruolo. Adoperiamo cioè la chiave privata per cifrare e quella pubblica per verificare. Il testo che viene scambiato tra Alice e Bob è quindi leggibile a chiunque, ma questo non rappresenta un problema perché l'obiettivo non è più la riservatezza sul contenuto del messaggio, ma l'assicurazione, per Bob, che il messaggio provenga proprio da Alice e che non sia stato modificato rispetto all'originale (autenticità ed integrità).

Anche per la *digital signature* i problemi derivanti dalla 'pesantezza' compu-tazionale della crittografia asimmetrica hanno reso necessario aggiungere dei passi in più, che, paradossalmente, rendono tutta l'operazione molto più 'leggera'. Introduciamo quindi un nuovo concetto: *l'hashing*

---

[30] In questo modo l'operazione di cifratura asimmetrica (3) è applicata alla sola sequenza, di lunghezza limitata, della chiave simmetrica; al giorno d'oggi stiamo parlando di stringhe (sequenze) lunghe al massimo 256 bit.

### 5.3. Funzioni di *hash*

Le funzioni di *hash* sono delle funzioni matematiche che calcolano per un testo di qualsiasi lunghezza una stringa di lunghezza fissa (detta *digest*).

La caratteristica peculiare delle funzioni di *hash* usate in crittografia, è di essere unidirezionali[31]. In altri termini, una volta calcolato il *digest* di un testo è praticamente impossibile risalire al testo originario. Analogamente, è praticamente impossibile produrre un testo che dia per risultato un predefinito *digest*.

Il *digest* finisce quindi per essere una sorta di impronta digitale del documento stesso.

Formalmente, una funzione di *hash* $H$ effettua una trasformazione di un messaggio $m$ di lunghezza arbitraria[32] che riceve in input generando come output un messaggio (*digest* o impronta) $h = H(m)$ di lunghezza fissa ridotta rispetto a quella di $m$

Le funzioni di *hash* utilizzate in crittografia soddisfano le seguenti proprietà:

1. il messaggio $m$ in ingresso può essere di qualsiasi lunghezza;
2. il messaggio $h$ in uscita ha sempre lunghezza fissa;
3. il calcolo di $h=H(m)$ è veloce e poco oneroso;
4. la trasformazione è monodirezionale (*one-way*).

Perché una 'buona' funzione di *hash* soddisfi queste proprietà, in particolare l'ultima, devono quindi valere certi requisiti. Se $n$ è la lunghezza in bit del risultato della funzione di *hash* la probabilità di collisione[33] della funzione $H(m)$ deve essere estremamente bassa. Non dovrebbe essere cioè possibile per un attaccante trovare una coppia di messaggi $m$ e $m'$ con $m \neq m'$ tale che $H(m) = H(m')$ impiegando un numero di tentativi minore di circa $2^{n/2}$.

Inoltre:

- dato un possibile valore di output $Y$, l'attaccante non dovrebbe essere in grado di trovare un input $x$ tale che $Y=H(x)$ impiegando meno di $2^n$ tentativi;

[31] Le funzioni di *hash* rientrano nella categoria delle funzioni *one-way*, ma senza *trapdoor*

[32] Quasi; gli algoritmi reali hanno dei limiti, però molto grandi rispetto alla lunghezza dei messaggi o dei documenti normalmente utilizzati.

[33] Per collisione si intende che due messaggi producono lo stesso *digest*. A tal proposito ricordiamo il cosiddetto paradosso del compleanno: dati $n$ interi presi a caso da una distribuzione che si suppone uniforme (il grado di uniformità è una misura della qualità di un algoritmo di *hash*) compresa tra 1 e $N$, la probabilità di trovarne due uguali approssima 1 per $n=N/2$.

- dato un messaggio *m*, l'attaccante non dovrebbe essere in grado di trovare un secondo messaggio *m'* tale che *H(m)=H(m')* impiegando meno di $2^n$ tentativi.

## 5.4. La *digital signature* in pratica

Dopo la parentesi necessaria per introdurre le funzioni di *hash*, possiamo vedere come viene effettivamente calcolata una *digital signature*.

Il mittente Alice:

1. calcola il *digest* $D=H(T)$ per il messaggio *T* da trasmettere;
2. cifra il *digest* *D* con la propria chiave privata producendo $D_c= Enc(D, A_{sk})$;
3. invia il testo *T* ed il *digest* $D_c$ cifrato.

Il destinatario Bob:

4. calcola il *digest* $D' = H(T')$ per il testo *T'* che ha ricevuto;
5. decifra il *digest* $D_c$ ricevuto utilizzando la chiave pubblica di Alice e ottenendo $D= Dec(D_c, A_{pk})$;
6. confronta *D* e *D'*.

Se i due *digest* sono uguali, il destinatario Bob ha la garanzia che il mittente è proprio Alice (da cui il nome 'firma digitale'). Ma non solo! Poiché anche la modifica di un solo bit del testo modificherebbe il valore del *digest*, in questo modo Bob ha anche la garanzia che *T=T'* cioè che il messaggio non è stato modificato nel tragitto percorso.

Per evidenziare quest'ultimo punto, mostriamo come due stringhe che differiscono esclusivamente di un carattere producano risultati completamente diversi quando sottoposte ad un trasformazione di *hash*.

| Testo[34] | *Digest* (SHA-1) |
|---|---|
| Italia-Germania 4-3 | 6CF5982ECB6BBE81882A03BC205D8857697C9AA6 |
| Italia-Germania 5-3 | A3BB835E86561F172E233A8F495E75E94FE640EE |

## 5.5. L'associazione tra la chiave pubblica ed il suo possessore

Un lettore accorto si sarà reso conto che nelle varie descrizioni abbiamo trascurato un paio di questioni peraltro essenziali.

Come fanno Alice e Bob a procurarsi le rispettive chiavi pubbliche?

Come fanno a sapere che quelle chiavi sono veramente di Alice e Bob?

Le soluzioni che si sono affermate sono sostanzialmente due: la PKI (Public Key Infrastructure, basato sulla figura di una *terza parte fidata* detta Certification Authority[35]) e il Pretty Good Privacy (PGP) e la sua architettura di web

[34] I due testi differiscono di un unico bit!

[35] E' il modello adottato dall'Italia per la firma digitale a valore legale.

trust.

Nel primo caso la *terza parte* emette un certificato digitale, che attesta che una chiave pubblica si riferisce ad un determinato soggetto (una persona, un router, un server, ecc.). Questo certificato, che potremmo paragonare ad un atto notorio davanti ad un Pubblico Ufficiale, ha un formato codificato da standard internazionali *de jure*[36] e *de facto*[37] ed è a sua volta firmato dalla *terza parte* a garanzia di integrità e autenticità. La chiave privata della Certification Authority diventa la chiave di volta di tutta la catena del trust e richiede quindi speciali misure di protezione.

Il concetto della Web of Trust è più semplice: uno che conosco e di cui mi fido (perché lavoriamo insieme, perché l'ho conosciuto durante una conferenza o per qualsiasi altro ragionevole motivo) e con cui ho scambiato fisicamente le chiavi pubbliche, mi assicura, firmandola, che la chiave che ho appena ricevuto per email è quella di Bob. Da quel momento, per la proprietà transitiva del Trust, io mi fido anche di Bob e di quelli che lui mi presenterà in futuro. Semplice e lineare, ma non facilmente scalabile, cioè portabile dai numeri piccoli a quelli che caratterizzano le grandi organizzazioni aziendali o governative.

Non possiamo addentrarci oltre nei dettagli perché la loro descrizione richiederebbe troppo spazio; consiglio al lettore interessato ad approfondire, di sfruttare le moltissime risorse di Internet, usando la wikipedia[38] come punto di partenza, o fare riferimento ai testi indicati in bibliografia ([1],[10]).

## 6. La robustezza degli algoritmi crittografici

Come abbiamo visto, alla base della crittografia, simmetrica o asimmetrica, sta la difficoltà di risolvere un problema di inversione di una funzione. Lo abbiamo analizzato con qualche dettaglio per la crittografia asimmetrica, visto nella definizione delle funzioni di *hash*, ma il concetto si applica tal quale anche alla crittografia simmetrica, in cui il *plaintext* è sottoposto ad una trasformazione dipendente dalla chiave usata, che deve poi essere invertita per riottenere il testo in chiaro dal *ciphertext*. La difficoltà di questa inversione è, grossolanamente, funzione della dimensione della chiave, che rappresenta la 'difficoltà' del problema.

Qual è quindi lo scopo del nostro attaccante che cerca di rubare preziose informazioni che transitano in rete o di decifrare i progetti di un nuovo dispositivo, memorizzati su qualche hard disk? Indovinare la chiave con cui questi dati

---

[36] *ISO/IEC 9594-8 - Information technology - Open Systems Interconnection - The Directory: Public-key and attribute certificate frameworks*, più noto come *ITU X.509 Information Technology – Open Systems Interconnection – The Directory: Authentication Framework; ITU-T Recommendation X.509* http://www.itu.int/rec/T-REC-X.509/en (visitato il 03/07/2010).

[37] Per chi non ha mai sentito parlare delle RFC (Request for comments), queste sono specifiche tecniche sviluppate dalla comunità internet: www.ietf.org.

[38] http://en.wikipedia.org/wiki/Public_key_infrastructure (visitato il 03/07/2010).

sono stati cifrati e portarli da *ciphertext* nuovamente a *plaintext*.

In passato si tendeva (e spesso si tende anche adesso, con risultati generalmente pessimi) a proteggere i dati con un algoritmo segreto. Sempre più, oggi si ritiene che la cosiddetta 'security by obscurity' aiuti sostanzialmente i "cattivi" che possono così più facilmente[39] ottenere il loro bottino, e che, al contrario, un algoritmo pubblico, le cui caratteristiche siano note a tutti e valutate dall'intera comunità mondiale dei crittografi, abbia una maggiore probabilità di non contenere *bugs* o trappole.

In conclusione, la capacità di un sistema crittografico di resistere ad un attacco dipende dalla 'difficoltà' del problema: la lunghezza della chiave.

Quindi, gli algoritmi crittografici, un po' come i cibi, scadono. La disponibilità di processori più potenti, di memoria a basso costo e di nuove idee per lo sfruttamento dei processori (per es. Grid computing[40]) rendono possibili, con tempi e costi ragionevoli,[41] calcoli inavvicinabili in precedenza. Da questo consegue che l'attacco concettualmente più semplice che si può portare all'algoritmo ('forza bruta', provare tutte le chiavi possibili) diventa realizzabile in concreto.

A ciò si aggiungono le astute scorciatoie che la ricerca matematica ha messo e mette continuamente a disposizione, e altre brillanti metodologie collaterali che sfruttano il fatto banale ricordato prima, che i computer sono oggetti fisici e che la manipolazione dei dati consuma energia e tempo.[42]

Analizzare questi aspetti ci porterebbe troppo lontano, ci limiteremo quindi a parlare brevemente della resistenza contro la forza bruta.

Se parliamo di chiavi simmetriche, il concetto dovrebbe essere abbastanza semplice. Una chiave è una stringa lunga $n$ bit, corrispondenti a $2^n$ possibili stringhe differenti. Se l'algoritmo di generazione delle chiavi (che non fa parte dell'algoritmo di cifratura) assicura una distribuzione uniforme (cioè se nel generare la chiave abbiamo a disposizione un buon generatore di numeri casuali, meglio se basato su una sorgente fisica realmente imprevedibile), la probabilità di azzeccare la chiave è $2^{-n}$ . Con la potenza di calcolo disponibile oggi si richiede che $n > 80$. Ma attenzione! Non dobbiamo dimenticarci di tener conto di che cosa stiamo proteggendo e per quanto tempo pensiamo debba rimanere segreto. Un segreto di Stato dovrà reggere almeno 30 anni, contro attaccanti determinati e probabilmente forniti di enorme potenza elaborativa, la password per accedere al Social Network deve resistere qualche mese e, in ogni caso, proteggerà, in genere, dati non troppo importanti.

Un po' più complicato è fare questa stessa valutazione per gli algoritmi asimmetrici. Le raccomandazioni [26] cambiano leggermente secondo i metodi adoperati che, in questo caso, tengono conto delle migliori tecniche matematiche disponibili per ciascun algoritmo.

[39] Che si adotti la sempreverde tecnica della corruzione o il più sofisticato *reverse engineering*, la segretezza dell'algoritmo non si riesce a preservare troppo a lungo; a quel punto è meglio che le sue debolezze eventuali siano note a tutti e non solo ai 'malvagi'.

[40] http://en.wikipedia.org/wiki/Grid_computing (visitato il 17/06/2010).

[41] Ragionevolezza che è funzione del tipo di dato che viene protetto e della capacità elaborativa presunta dell'attaccante.

[42] http://en.wikipedia.org/wiki/Side_channel_attack (visitato il 17/06/2010).

Per fare un esempio, la società RSA, che prende il nome dai tre inventori dell'omonimo algoritmo e che ne possedeva il brevetto (ormai scaduto), dal 1999 al 2007 ha lanciato una serie di sfide, fornendo numeri di lunghezza crescente da fattorizzare. Al momento il record è stata la fattorizzazione di un numero lungo 768 bit.[43] Le chiavi adoperate per la firma digitale usano chiavi da 1024 bit, hanno quindi ancora qualche anno davanti di tranquillità, anche se già oggi ci stiamo predisponendo per andare oltre.

Ad oggi la principale debolezza del sistema della firma è emersa dall'algoritmo di *hash* più diffuso[44] (SHA-1), debolezza ancora lontana dal mettere a rischio i sistemi che su di esso si basano, ma che ha innescato il processo di migrazione verso algoritmi che producono digest più lunghi e, in una prospettiva più a lungo termine, verso una nuova famiglia di funzioni di *hash*.[45]

## 7. Le evoluzioni

### 7.1 Le curve ellittiche

La matematica delle curve ellittiche è decisamente più sofisticata di quella usata per la fattorizzazione, quindi ne daremo solo un breve accenno, indicando però i motivi per cui suscitano un grande interesse per le applicazioni pratiche.

Una curva ellittica è un'equazione della forma

$$y^2 = x^3 + ax + b$$

con $a$ e $b$ tali che $4a^3 + 27b^2 \neq 0 \ (mod \ p)$

Nel nostro contesto crittografico tutte le operazioni sono eseguite modulo $p$ (primo dispari[46]).

Notiamo che la curva è simmetrica rispetto all'asse delle $x$ (l'equazione è quadratica in $y$).

[43] Per maggiori dettagli http://www.rsa.com/rsalabs/node.asp?id=2092 (visitato il 18/06/2010).

[44] Per maggiori dettagli ed i link alle specifiche più tecniche http://it.wikipedia.org/wiki/Secure_Hash_Algorithm (visitato il 18/06/2010).

[45] La selezione del nuovo algoritmo *hash* di riferimento è gestita dal già citato NIST, con un concorso pubblico: http://csrc.nist.gov/groups/ST/hash/sha-3/index.html (18/06/2010).

[46] C'è anche un'altra tipologia, con i calcoli eseguiti nel *campo* finito GF($2^m$), $m$ intero

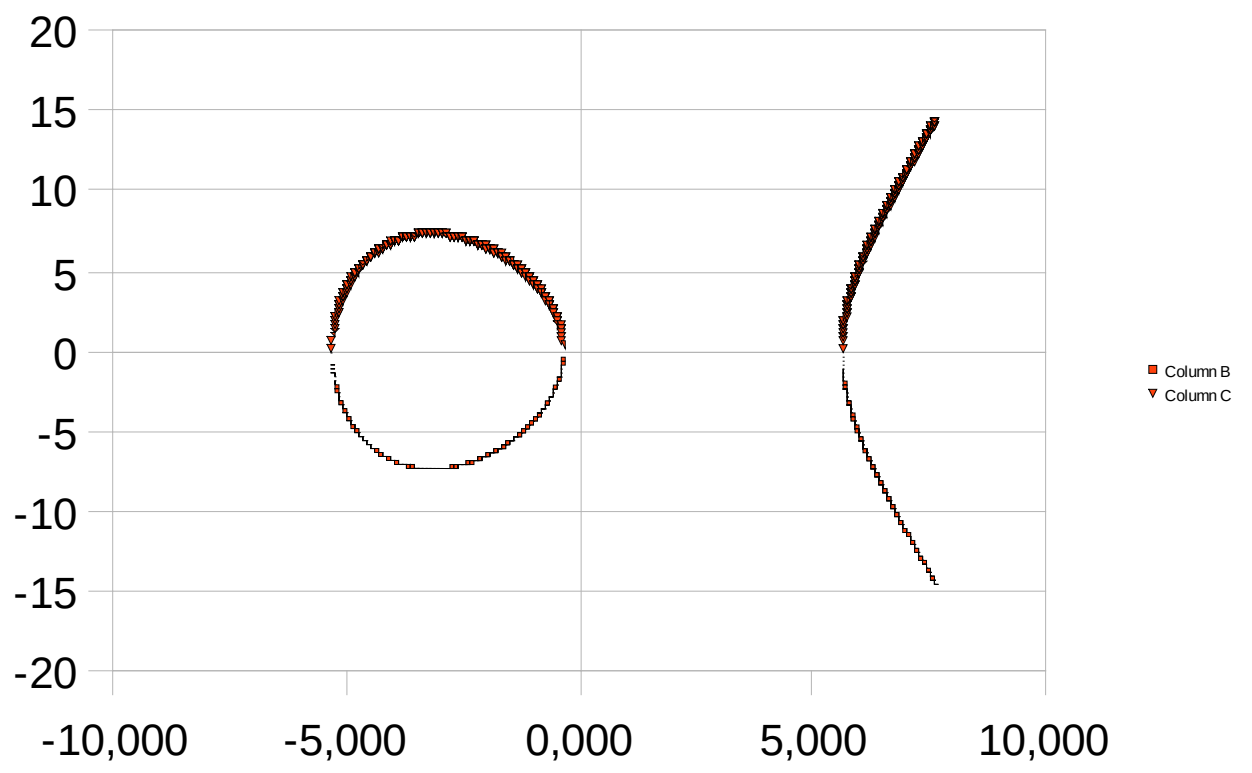

Figura 2 – Un esempio di curva ellittica (*a=-30* e *b=-10)*

Presi due punti sulla curva ($P$ e $Q$), il segmento che li unisce incontra la curva in un terzo punto (poiché la curva è di terzo grado). Si prende quindi il punto simmetrico rispetto all'asse $x$ di questo terzo punto e si definisce il nuovo punto $R=P+Q$.

Si trova che l'addizione così definita, con qualche altro accorgimento, soddisfa tutte le proprietà algebriche usuali dell'addizione.

Scegliendo $P$ e $Q$ coincidenti si definisce $P+P=2P$ e per estensione $kP$ e $P^n$

### 7.1.1 - Qual'è il problema da risolvere

Il problema matematico che deve essere risolto è, dati due punti $Q$ e $P$, tali che $Q=kP$, trovare il valore di $k$.

Un possibile approccio è quello, banale, di andare per tentativi provando tutti i valori di $k$ fino a trovare quello giusto. Risulta che, se il valore del primo $p$ è sufficientemente grande, il numero di tentativi può diventare enorme.

Esistono delle tecniche che permettono di ridurre il numero di tentativi, che tuttavia rimane grande ( $O\left(\sqrt{p}\right)$.

### 7.1.2 - In quali contesti la crittografia a curve ellittiche può essere utile

Il principale vantaggio offerto dall'uso delle curve ellittiche sta nella ridotta lunghezza delle chiavi a parità di sicurezza. Una chiave ECC[47] a 160 bit offre la stessa sicurezza di una RSA a 1024 bit, occupando meno memoria e richiedendo meno energia per le operazioni di firma, il che rende questa crittografia adatta all'uso in dispositivi portatili (cellulari e palmari) dove ci sono costrizioni derivanti dalla capacità elaborativa del processore, dalle batterie e dalla memoria disponibile.

## 7.2 - La crittografia quantistica

[47] Elliptic Curve Cryptography

La crittografia quantistica[48] è una forma alternativa di distribuzione delle chiavi, basata sui principi alla base della fisica quantistica, cioè sul fatto che, semplicemente effettuando una misura, il sistema misurato ne viene perturbato.

Per descrivere la crittografia quantistica, occorrerebbe richiamare alcuni concetti della fisica sottostante e si andrebbe ben oltre lo scopo di queste note.

Rimanendo a livello generale, lo scambio di chiavi quantistico permette di rendere trascurabile la probabilità che Eve possa avere accesso alla chiave senza che Alice e Bob non se ne accorgano. Grazie a questa proprietà è possibile, in teoria, applicare la già citata cifratura *One-Time pad*.

Peraltro, di recente è stato provato con successo un attacco ad un sistema già operativo di crittografia quantistica[49], sfruttando le inevitabili differenze tra i modelli matematici e la loro realizzazioni effettiva.

Anche senza considerare quest'ultima notizia, è bene sottolineare che se la crittografia vera e propria che poi viene adoperata negli scambi informativi è la crittografia classica simmetrica, la sicurezza complessiva del sistema viene comunque a dipendere dalla robustezza degli algoritmi crittografici[50].

E, *last but not least*, come per tutte le applicazioni informatiche, vanno considerate l'effettiva traduzione dell'algoritmo in un programma eseguibile e la sicurezza delle stazioni di lavoro su cui i dati vengono elaborati.

## 8.3 - I computer quantistici

Il computer quantistico è un qualcosa ancora inesistente, ma che se e quando sarà effettivo, cambierà drasticamente il panorama della crittografia e dell'informatica in generale.

Infatti, oltre a permettere un ulteriore balzo verso la miniaturizzazione dei circuiti, gli stessi principi fisici che governano il mondo atomico e subatomico[51] permettono in teoria una modalità di eseguire i calcoli in modo molto più efficiente.

Anche in questa sezione, la descrizione non potrà che essere sommaria poiché la fisica sottostante non è semplicissima e, purtroppo, non fa parte del programma delle scuole superiori, nemmeno a livello puramente descrittivo.

Il principio fondamentale della fisica quantistica è che un sistema può esiste-re in una sovrapposizione di stati. Quindi, mentre un bit tradizionale può assumere i valori 0 e 1, un qubit, come viene denominata l'unità di informazione di un computer quantistico, può essere nello stato |0> oppure |1> oppure in una sovrapposizione dei due. Da questa proprietà ne deriva che un insieme di $N$ qubit può assumere $2^N$ stati differenti allo stesso tempo e il calcolo quantistico diventa analogo ad un calcolo effettuato su un computer tradizionale da $2^N$ processori che lavorano in parallelo.

---

48 http://arxiv.org/PS_cache/quant-ph/pdf/0101/0101098v2.pdf (02/07/2010)

49 http://www.nature.com/news/2010/100520/full/news.2010.256.html (2/7/2010)

50 Se si decide di usare un algoritmo standard come l'AES, la sicurezza complessiva è data dalla robustezza di quest'ultimo, non importa quanto sicuro sia stato il processo di scambio delle chiavi

51 In realtà la fisica quantistica dovrebbe applicarsi a tutti gli oggetti, ma nei sistemi macroscopici gli effetti quantistici diventano trascurabili. Si noti che questa affermazione sintetizza in maniera assai grossolana un dibattito ancora aperto sull'interpretazione della meccanica quantistica, su alcune debolezze logiche che ancora non sono state chiarite, a dispetto dell'enorme mole di dati che ne confermano la validità. Da alcune risposte "interpretative" potrebbe anche derivare la fattibilità o meno di un computer quantistico.

Da questa caratteristica deriva, per esempio, che il tempo di fattorizzazione di un numero primo si riduce di parecchi ordini di grandezza[52]. Tuttavia non tutti i problemi "complessi" risulteranno semplificati dal computer quantistico[53].

E, inoltre, non è affatto scontato che il computer quantistico risulti effettivamente realizzabile. I problemi legati a mantenere un sistema macroscopico composto da migliaia di qubit in uno stato quantistico è di una difficoltà enorme. Molti progressi sono stati fatti, ma la strada da fare è ancora lunga.

## Appendice 1 - Un po' di matematica (quasi) elementare

Per comprendere i concetti base della crittografia asimmetrica, è necessario introdurre qualche concetto che non è scontato faccia parte del normale curriculum scolastico.

Come vedremo, per afferrare le idee principali della crittografia a chiave pubblica sono sufficienti le nozioni imparate alle medie, con qualche escursione nell'algebra delle superiori per la dimostrazione di qualche teorema.

Infatti, uno degli aspetti più affascinanti della Teoria dei Numeri è che capirne i problemi (gli enunciati) è alla portata di tutti, laddove risolverli è, in molti casi, di una difficoltà e complessità incredibile. Pensate al Teorema di Fermat[54] o alla congettura di Goldbach:[55] il primo ha dovuto attendere oltre tre secoli per essere dimostrato e tecniche assai sofisticate; la seconda, verificata fino a $3 \cdot 10^{17}$ è ritenuta vera, ma non è stata (ancora?) dimostrata e c'è anche chi ritiene che non sia dimostrabile.[56]

Alla base di tutto ci sono gli 'atomi' dell'aritmetica, gli elementi più semplici

52 http://en.wikipedia.org/wiki/Shor%27s_algorithm (21/5/2010)

53 S. Aaronson, I limiti del computer quantistico, Le Scienze, maggio 2008

54 Dati $a, b, c, n \in \mathbb{N}$ la relazione $a^n + b^n = c^n$ non ha soluzioni per n>2 ; Per n=2, ritroviamo il Teorema di Pitagora.

55 Ogni numero pari > 2 si può esprimere come la somma di due numeri primi.

56 L'indimostrabilità potrebbe derivare da un altro celeberrimo teorema, noto anche fuori dall'ambito dei matematici, quello di Gödel sull'incompletezza dell'aritmetica: in un sistema formale (come l'aritmetica), esistono proposizioni che non sono derivabili o dimostrabili nel sistema stesso, pur essendo "vere" (D. Knuth *All Questions Answered* http://www.ams.org/notices/200203/fea-knuth.pdf. pag. 4, visitato il 21/07/2010).

che stanno alla base dei concetti che verranno esposti: i numeri interi[57] e i numeri primi.

Come è noto, indicando con $p_i$ un generico numero primo, ogni numero naturale $N \in \mathbb{N}$ può essere espresso:

$$N = p_1^{n_1} \cdot p_2^{n_2} \cdot p_3^{n_3} \cdot \ldots p_m^{n_m} \quad (10)$$

con <u>una</u> sola combinazione di primi ed esponenti.

I numeri primi sono infiniti[58].

Questo risultato fu già dimostrato da Euclide, con una prova tanto semplice quanto elegante (vedi riquadro).

**La dimostrazione dell'infinità dei numeri primi**

Per iniziare notiamo che se $d$ divide due interi $a$ e $b$ , allora divide anche la differenza ($a - b$). Ora supponiamo che i primi siano un numero finito $n$, $p_1 < p_2 < p_3 < \cdots < p_n$ , $p_n$ il più grande.

Prendiamo i due numeri:

$$a = p_1 \cdot p_2 \cdot p_3 \ldots \cdot p_n \quad \text{e} \quad b = p_1 \cdot p_2 \cdot p_3 \ldots \cdot p_n + 1$$

$b$ non può essere primo perché, per ipotesi, $p_n$ è il più grande numero primo, quindi deve essere divisibile per qualche $p_i$, $a$ è ovviamente divisibile per lo stesso $p_i$ (per costruzione), quindi $(b - a) = 1$ è divisibile per $p_i$; questo è impossibile quindi l'ipotesi iniziale (la finitezza di $n$) è falsa.

Assodato che di numeri primi ne abbiamo a sufficienza, la domanda successiva è 'Come sono distribuiti?'. A questo risponde il Teorema dei Numeri Primi[59], la cui dimostrazione è invece piuttosto tecnica e non daremo. Se indichiamo con $\pi(x)$ il numero di primi minori o uguali ad un numero reale $x$, il teorema asserisce che per $x$ molto grande (*tende all'infinito*) si ha che:

.

In maniera più precisa:

$$\lim_{x \to \infty} \frac{\pi(x)}{x / \ln(x)} = 1 \quad (11)$$

Per capire l'importanza pratica di questo teorema, proviamo a calcolare quanti primi ci sono tra $2^{1023}$ e $2^{1024}$, cioè nell'intervallo di numeri rappresentato da una

57 "Dio ha creato i numeri interi, tutto il resto è opera dell'uomo", L. Kronecker.

58 Il che ci rassicura che ce ne siano a sufficienza per le esigenze della nostra crittografia.

59 Il teorema ha una lunga storia e molti padri: Gauss e Legendre, indipendentemente, lo espressero come congettura senza dimostrarlo; la prima dimostrazione risale al 1896 (Hadamard e de la Vallée-Poussin, di nuovo indipendentemente, usando i metodi dell'analisi complessa) fino a Erdos e Selberg, che lo dimostrarono usando metodi 'elementari'.

sequenza di 1024 bit, la ‘lunghezza’ attualmente adottata per le chiavi della firma digitale.

I numeri in questione sono sufficientemente grandi da assumere la relazione quasi esatta; quindi, si ottiene che nell'intervallo indicato ci sono all'incirca $10^{305}$ numeri primi.

Per darvi un'idea, $10^{80}$ è il numero stimato di atomi dell'universo!!

**La fattorizzazione è un problema ‘difficile’**

Ognuno di voi ricorderà certamente il concetto di minimo comune multiplo, utilizzato a scuola essenzialmente nel calcolo frazionario, per dare alle frazioni lo stesso denominatore (minimo comune denominatore) e quindi agire esclusivamente sui numeratori per sommarle o confrontarle.

Per fare un esempio veramente banale:

E' vera la seguente ipotesi?

$$\frac{3}{7} > \frac{4}{9}$$

Per verificarla si trova il minimo comune multiplo (mcm) dei denominatori.

Vi ricordo che per calcolare il minimo comune multiplo di un set di numeri (*a,b,c,d,...*) si effettua la scomposizione[60] in fattori primi di ciascun numero:

$$a = p_1^n \cdot p_2^m \cdot ... p_n^{q...}$$

(dove i $p_i$ sono tutti primi)

<u>Il mcm è il numero che si ottiene moltiplicando tutti i fattori primi dei numeri dati, presi una sola volta e con il massimo esponente.</u>

Tornando al nostro esempio, 7 è un numero primo, $9=3^2$, quindi:

$$mcm = 7 \cdot 3^2 = 63$$

Le nostre frazioni diventano quindi

$$\frac{27}{63} > \frac{28}{63}$$

quindi l'ipotesi era falsa

Chiaro e semplice.

Ma che cosa accade se il numero di cui dobbiamo trovare la scomposizione in fattori primi è 171371?

Non è pari e questo esclude la divisibilità per 2. Non finisce per 5 o 0, quindi non è divisibile per 5. La somma delle cifre non è un multiplo di 3 quindi non è divisibile per 3. La somma delle sue cifre in posizione pari meno la somma delle cifre in posizione dispari non è zero né un multiplo di 11, quindi non è divisibile per 11.

Diciamo che qui si esauriscono i criteri ‘semplici’ di divisibilità. Dopo, si procede essenzialmente per tentativi, sia pure con metodi sofisticati che evitano di dover provare tutte le possibilità, metodi che non modificano però la com-

60 Abbiamo visto che questa scomposizione è unica plessità del problema.
Torniamo al nostro numero. Per chi volesse provarci prima di passare alla soluzione[61] di seguito sono elencati i numeri primi minori di 1000.

*Tabella 2*

```
For more information on primes see http://primes.utm.edu/

2 3 5 7 11 13 17 19 23 29 31 37 41 43 47 53 59 61 67 71 73
79 83 89 97 101 103 107 109 113 127 131 137 139 149 151
157 163 167 173  179 181 191 193 197 199 211 223 227 229
233 239 241 251 257 263 269 271 277 281 283 293 307 311
313 317 331 337 347 349 353 359 367 373 379 383 389 397
401 409 419 421 431 433 439 443 449 457 461 463 467 479
487 491 499 503 509 521 523 541 547 557 563 569 571 577
587 593 599 601 607 613 617 619 631 641 643 647 653 659
661 673 677 683 691 701 709 719 727 733 739 743 751 757
761 769 773 787 797 809 811 821 823 827 829 839 853 857
859 863 877 881 883 887 907 911 919 929 937 941 947 953
967 971 977 983 991 997
```

Il risultato cercato è $409 \cdot 419$, ma chi ci ha provato si sarà reso conto della grande differenza che passa tra moltiplicare questi due numeri e individuarli per tentativi dal loro prodotto.

Al momento l'unico modo conosciuto per trovare i fattori primi di un numero è, sostanzialmente, quello di andare per tentativi, con un tempo di elaborazione che cresce con $n$ in maniera sub esponenziale.[62]

Tornando alla nostra storia, è immaginabile la complessità di doverlo fare per un numero lungo 300 cifre decimali (a questo corrisponde un numero binario di 1024 bit), dovendo inoltre determinare se i numeri che adoperate per provare siano effettivamente primi.[63]

## L'aritmetica modulare (congruenze)

Prima che vi spaventiate, sappiate che si tratta di una modalità di fare le somme che vi è assolutamente consueta, addirittura quotidiana.

---

[61] Suggerimento: i numeri da provare sono quelli minori della $\sqrt{\ }$ del numero esaminato.

[62] Cioè più lentamente di $a^n$ per qualche $a > 1$, ma più velocemente di qualsiasi polinomio in $n$.

[63] Solo di recente il problema della 'primalità' di un numero è stato risolto, dimostrando che si tratta di un problema 'abbordabile' cioè che è un problema di classe P ed è stato descritto l'algoritmo che ne permette la risoluzione. Anche così, comunque, la scomposizione in fattori primi di un numero grande rimane un problema difficile.

Per capirci, sono le 10 di sera e state impostando la sveglia, siete stanchi e volete dormire almeno 8 ore. A che ora la puntate? Le 10 di sera corrispondono alle 22, quindi

$$22+8=30$$

e poiché il giorno è di 24 ore, la sveglia viene puntata alle 6 (30-24).

Bene, questo è un esempio di somma modulo 24 e si esprime in questa forma:

$$22+8\equiv 6\ (mod\ 24)$$

In generale, dire

$$a\equiv b\ (mod\ n) \qquad (12)$$

(si legge a congruente a b modulo n) vuol dire che

$$a=k\cdot n+b \qquad (13)$$

con *k* un qualsiasi numero intero. Detto ancora in altre parole, *b* è il resto della divisione *a:n.*

Il motivo per introdurre l'aritmetica modulare è che, in questa maniera, i numeri da trattare restano di dimensione limitata. Quando si considerano i numeri in gioco nell'algoritmo RSA, l'uso del modulo permette di avere a che fare con numeri di dimensione controllata ed è anche necessario per la robustezza matematica dell'algoritmo stesso.

L'aritmetica modulare, applicata ai numeri interi relativi, ha proprietà quasi uguali a quelle dell'aritmetica tradizionale. Ritroviamo le proprietà commutativa, associativa e distributiva dell'addizione e della moltiplicazione.

proprietà commutativa dell'addizione e della moltiplicazione

$$a+b\equiv c\ (mod\ m)\Rightarrow b+a\equiv c\ (mod\ m) \qquad (14)$$

$$a\cdot b\equiv c\ (mod\ m)\Rightarrow b\cdot a\equiv c\ (mod\ m) \qquad (15)$$

proprietà associativa dell'addizione e della moltiplicazione

$$((a+b)+c)\equiv d\ (mod\ m)\Rightarrow (a+(b+c))\equiv d\ (mod\ m) \qquad (16)$$

$$((a\cdot b)\cdot c)\equiv d\ (mod\ m)\quad\Rightarrow\quad (a\cdot(b\cdot c))\equiv d\ (mod\ m) \qquad (17)$$

proprietà distributiva

$$a\cdot(b+c)\equiv d\ (mod\ m)\quad\Rightarrow\quad a\cdot b+a\cdot c\equiv d\ (mod\ m) \qquad (18)$$

La divisione, contrariamente all'aritmetica tradizionale, non fornisce un risultato frazionario, ma restituisce sempre un intero compreso tra 0 e *m*-1, dove *m* è il modulo.

Ricordiamo infatti che la divisione è l'operazione inversa della moltiplicazione e quindi

$$\frac{a}{b} \equiv d \,(mod\ m) \quad \Rightarrow \quad a \equiv b \cdot d \,(mod\ m) \tag{19}$$

Esempio:

Quanto vale $\frac{4}{7}\,(mod\,11)$? E' facile verificare, provando i vari valori, che il risultato è 10. Infatti $7 \cdot 10 = 70$ che equivale a $6 \cdot 11 + 4$.

Troviamo ancora che valgono altre proprietà quasi ovvie:

1) l'esistenza di un elemento neutro, cioè di un elemento per il quale l'operazione non produce cambiamenti:

$$(a+0)(mod\ m) \equiv (0+a)\ (mod\ m) \equiv a\,(mod\ m). \tag{20}$$

$$(a \cdot 1)(mod\ m) \equiv (1 \cdot a)(mod\ m) \equiv a\,(mod\ m). \tag{21}$$

2) l'esistenza di un inverso, cioè che dati due numeri interi qualsiasi *a* e *m*, esiste un numero intero (-*a*) tale che sommato ad *a* dà come risultato l'elemento neutro dell'addizione:

$$(a + (-a))\ (mod\ m) \equiv 0\,(mod\ m) \tag{22}$$

Analogamente, dati due numeri interi qualsiasi *a* e *m* primi tra loro[64] con $a \neq 0$, esiste un numero intero $a^{-1}$ tale che moltiplicato per *a* dà come risultato l'elemento neutro della moltiplicazione:

$$(a \cdot a^{-1})\ (mod\ m) \equiv 1\ (mod\ m) \tag{23}$$

3) La proprietà di chiusura rispetto alle operazioni di addizione e moltiplicazione

$$\begin{aligned} a,b \in G \quad &\Rightarrow \quad (a+b) \in G \\ a,b \in G \quad &\Rightarrow \quad (a \cdot b) \in G \end{aligned} \tag{24}$$

4) La proprietà transitiva

[64]quando MCD ($m_i$, $m_j$)=1, $m_i$, $m_j$ si dicono relativamente primi; qui e nel seguito MCD sta per Massimo Comun Divisore e mcm per minimo comune multiplo.

$$a \equiv b \ (mod\ m) \quad e \quad a \equiv c \ (mod\ m) \quad \Rightarrow \quad b \equiv c \ (mod\ m) \qquad (25)$$

Una delle più importanti caratteristiche non ovvie del calcolo modulare è la riducibilità degli addendi e dei fattori:

$$a + b \equiv c \ (mod\ m) \quad \Rightarrow \quad (a \ (mod\ m) + b \ (mod\ m)) \equiv c \ (mod\ m) \qquad (26)$$

$$a \cdot b \equiv c \ (mod\ m) \quad \Rightarrow \quad (a \ (mod\ m) \cdot b \ (mod\ m)) \equiv c \ (mod\ m) \qquad (27)$$

Ad esempio, supponiamo di voler calcolare $2^{10}$ mod 11. Possiamo scrivere $2^{10} = 2^4 \cdot 2^4 \cdot 2^2$.
Possiamo quindi ridurre

$$2^4 = 16 \equiv 5 \ (mod\ 11)$$
$$2^2 = 4 \equiv 4 \ (mod\ 11)$$

e la nostra espressione diventa

$$5 \cdot 5 \cdot 4 = 25 \cdot 4 \equiv 3 \cdot 4 \equiv 1 \ (mod\ 11) \qquad (28)$$

Per verificare: $2^{10} = 1024 = 11 \cdot 93 + 1$

**Ulteriori proprietà delle congruenze**

Congruenze con lo stesso modulo possono essere aggiunte tra loro e moltiplicate, come nell'aritmetica tradizionale, ma con qualche sottile differenza:

$$a \equiv b \ (mod\ m) \text{ e } c \equiv d \ (mod\ m) \Rightarrow (a + c) \equiv (b + d)(mod\ m)$$
$$\Rightarrow (a \cdot c) \equiv (b \cdot d)(mod\ m) \qquad (29)$$
$$\Rightarrow a^n \equiv b^n \ (mod\ m)$$

Se $n = k \cdot d$ allora $a \equiv b \ (mod\ n) \Rightarrow a \equiv b \ (mod\ d)$
Sia M=mcm($m_1$,$m_2$); si ha

$$a \equiv b \ (mod\ m_1) \text{ e } a \equiv b \ (mod\ m_2) \Leftrightarrow a \equiv b \ (mod\ M)$$

Più in generale,

$$a \equiv b \ (mod\ m_i) \text{ e } \forall\, i,j \in [1, \ldots k] \ i \neq j \ MCD\ (m_i, m_j) = 1$$
$$\Leftrightarrow a \equiv b \ (mod\ m_i m_j \ldots m_k) \qquad (30)$$

Poiché ogni $m$ composto può essere scritto come $m = p_1^{e_1} \cdot p_2^{e_2} \ldots p_k^{e_k}$ si ha in-

fine che

$$a \equiv b \ (mod \ m) \Leftrightarrow a \equiv b \left( mod \ p_i^{e_i} \right) \quad \forall i \tag{31}$$

Se $a \cdot k \equiv b \cdot k \ (mod \ n)$ il fattore comune $k$ non può essere semplicemente rimosso, come faremmo nell'aritmetica ordinaria, se non quando $k$ e $n$ sono primi tra loro, cioè quando MCD ($k$,$n$) = 1.

In generale, se MCD ($k$,$n$) = $d$, si ha

$$a \cdot k \equiv b \cdot k (mod \ n) \Rightarrow a \equiv b \left( mod \ \frac{n}{d} \right) \tag{32}$$

Infine, se $a$, $b$, $n$ sono tutti divisibili per qualche $d$, allora si ha

$$a \equiv b \ (mod \ n) \Rightarrow \frac{a}{d} \equiv \frac{b}{d} \left( mod \ \frac{n}{d} \right) \tag{33}$$

**La funzione φ di Eulero**

Per ogni intero $m$, scegliamo gli $n_i \in \{1, 2, ..., m-1\}$ che sono relativamente primi rispetto a $m$. La funzione $\varphi(m)$ di Eulero conta quanti sono questi numeri. Vale la pena di sottolineare che se $m$ è primo, tutti i numeri inferiori a $m$ sono primi rispetto a $m$ e quindi

$$\forall \ m \ primi \quad \varphi(m) = m - 1 \tag{34}$$

**Proprietà dellaφ di Eulero**

Se $m = p^{\alpha}$ , gli unici numeri che hanno un fattore comune con $m$ sono

$$p, 2p, \ldots, p^{\alpha-1} p \tag{35}$$

i quali sono esattamente $p^{\alpha-1}$

Quindi

$$\varphi\left(p^{\alpha}\right) = p^{\alpha} - p^{\alpha-1} = p^{\alpha} \cdot \left(1 - \frac{1}{p}\right) \tag{36}$$

In maniera appena un po' più complicata, ma seguendo la stessa falsariga, si

dimostra che se

$$m = p_1^{e_1} \cdot p_2^{e_2} \dots p_n^{e_n} \qquad (37)$$

$$\varphi(m) = p_1^{e_1} \cdot \left(1 - \frac{1}{p_1}\right) \cdot p_2^{e_2} \cdot \left(1 - \frac{1}{p_2}\right) \dots \cdot p_n^{e_n} \cdot \left(1 - \frac{1}{p_n}\right) = \varphi\left(p_1^{e_1}\right) \cdot \varphi\left(p_2^{e_2}\right) \dots \cdot \varphi\left(p_n^{e_n}\right) \qquad (38)$$

Con questi elementi siamo ora in grado di introdurre il teorema di Eulero-Fermat e poi la riduzione degli esponenti nelle congruenze, che completano l'armamentario per capire l'algoritmo RSA.

**Il teorema di Eulero-Fermat**

Il teorema di Eulero-Fermat[65] afferma che presi due numeri naturali $a$ e $m$, primi tra loro, vale la relazione

$$a^{\varphi(m)} \equiv 1 \left(mod\ m\right) \qquad (39)$$

La dimostrazione è abbastanza semplice.

Siano

$$r_i \in \{1, 2, \dots, m-1\} \qquad (40)$$

i numeri [minori di] e [relativamente primi a] $m$ con $i = 1, 2, \dots, \varphi(m)$.

Prendiamo un numero $a$ qualsiasi, anch'esso relativamente primo a $m$, moltiplichiamolo per tutti gli $r_i$ e dividiamo ciascun prodotto per $m$.

Si ottiene

$$a \cdot r_i = q_i \cdot m + r'_i \qquad (41)$$

con $r'_i$ è relativamente primo a $m$, altrimenti il termine a destra della (36) sarebbe divisibile per $m$, questo implicherebbe che anche il termine a sinistra è divisibile per $m$, il che non è possibile per costruzione.

Quindi $r'_i \in \{1, 2, \dots, m-1\}$. Inoltre gli $r'_i$ sono tutti diversi tra loro.

Per verificarlo, riscriviamo la (41) in forma di congruenza:

$$a \cdot r_i \equiv r'_i \quad \left(mod\ m\right) \qquad (42)$$

---

[65] Il sottocaso enunciato da Fermat è anche chiamato il ‘piccolo teorema di Fermat’ per distinguerlo dal più famoso ‘Ultimo Teorema’. Come nel caso più famoso, Fermat non ne diede una dimostrazione. La dimostrazione di Eulero si applica al caso più generale.

e notiamo che, per definizione, $\forall (i, j)(i \neq j)\, r_i \neq r_j$ .

Se si avesse per qualche coppia $(i, j)$ $r'_i = r'_j$ $(i \neq j)$ allora $r'_i \equiv r'_j \pmod m \;\Rightarrow\; a \cdot r_i \equiv a \cdot r_j \pmod m$ ed essendo $a$ relativamente primo a $m$, si può cancellare e si avrebbe $r_i = r_j$ , cadendo in contraddizione.

Quindi sia gli $r_i$ che gli $r'_i$ coprono l'intero insieme dei [numeri minori di] e [relativamente primi a] $m$.

Moltiplicando tra loro tutte le congruenze (42), si ha

$$a^{\varphi(m)} \cdot r_1 \cdot r_2 \cdot \ldots r_{\varphi(m)} \equiv r'_1 \cdot r'_2 \cdot \ldots r'_{\varphi(m)} \pmod m \qquad (43)$$

Poiché il prodotto degli $r_i$ e degli $r'_i$ è lo stesso ed è (per costruzione) relativamente primo a $m$, si può elidere e si ottiene la (39).

Nel caso particolare che $m$ sia primo, $\varphi(m) = m - 1$ e abbiamo $a^{m-1} \equiv 1 \pmod m$ che è la forma con cui lo enunciò Fermat.

E' importante notare che questa relazione è vera per qualsiasi primo, ma è vera anche per molti numeri composti e non può essere adoperata per determinare se un numero è primo o meno.

## La riduzione delle potenze, le radici e i logaritmi

La possibilità di ridurre un'espressione, rendendola più semplice da manipolare, esiste anche per l'elevamento a potenza, con regole leggermente differenti:

$$a^b \pmod m \equiv a^{b \,(mod\ \varphi(m))} \pmod m \qquad (44)$$

La (44) si ottiene facilmente, esprimendo $b$ nella forma:

$$b = k \cdot \varphi(m) + q \qquad (45)$$

da cui

$$a^b = a^{k\varphi(m)+q} = \underbrace{a^{\varphi(m)} \cdots a^{\varphi(m)}}_{k\ volte} a^q \qquad (46)$$

Sfruttando la (27) e la (39), si riottiene la (44)

La definizione dell'estrazione di una radice o del logaritmo è identica (o quasi) a quella dell'aritmetica tradizionale

Il problema si può esprimere nella stessa maniera:

$$a \equiv b^c \pmod m \qquad (47)$$

Il problema del logaritmo è dato dal trovare $c$, mentre il problema della radice è dato dal trovare $b$, essendo noti gli altri due parametri.

## Appendice 2 - La sicurezza delle informazioni: i concetti base

Diamo ora brevemente qualche definizione per fissare di cosa parliamo quando parliamo di sicurezza delle informazioni, il campo di applicazione delle tecniche crittografiche.

Alcuni di questi concetti hanno una storia millenaria, altri sono di origine piu recente, divenuti importanti nell'era di Internet, quando l'IT (Information Technology) è diventata ICT (Information and Communication Technology) e la trasmissione dei dati è diventata una parte fondamentale della vita dei cittadini e delle imprese.

**Integrità**: la riduzione a livelli accettabili del rischio che possano avvenire cancellazioni o modifiche di informazioni a seguito di interventi di entità non autorizzate o del verificarsi di fenomeni non controllabili (come il deteriorarsi dei supporti di memorizzazione, degradazione dei canali trasmissivi, guasti agli apparati, problemi di approvvigionamento energetico, incendi, allagamenti, etc.).

**Riservatezza**: La riduzione a livelli accettabili del rischio che un'entità possa, volontariamente o involontariamente, accedere all'informazione senza esserne autorizzata.

**Disponibilità**: la riduzione a livelli accettabili del rischio che possa essere impedito alle entità autorizzate l'accesso alle informazioni a seguito di interventi di altre entità non autorizzate o del verificarsi di fenomeni non controllabili.

**Autenticazione**: La riduzione a livelli accettabili del rischio che l'interlocutore di uno scambio di dati non sia effettivamente colui che dichiara di essere.

**Non ripudio**: La riduzione a livelli accettabili della possibilità che una delle entità coinvolte in uno scambio di comunicazioni possa, in seguito, negare di aver partecipato allo scambio e di aver trasmesso certe informazioni.

Avrete certamente notato l'ossessiva ripetizione dell'espressione 'riduzione a livelli accettabili'. E' fondamentale che ognuno di noi si renda conto che la sicurezza assoluta, al 100%, non esiste, nell'ICT e nella vita di tutti i giorni. Per quante precauzioni si possano prendere, per quante contromisure si possano adottare, l'evento inatteso, dannoso, non voluto potrà sempre accadere. Quello che si può fare, nella vita e nel lavoro, è cercare di ridurre la probabilità che accada e di ridurne le conseguenze, spendendo tempo e denaro in maniera ragionevole.

Per fare un esempio, preso dalla vita di ogni giorno: esiste una probabilità non nulla che l'auto su cui viaggiamo per andare al lavoro si guasti; possiamo fare la manutenzione periodica, tenerla in garage, tuttavia un guasto o una foratura possono accadere, magari proprio il giorno in cui abbiamo fretta. Potremmo come soluzione comprare più macchine e tenerle parcheggiate lungo il percorso, una ogni 500 metri. Un costo niente male. Poi dovremmo essere certi che funzionino, quindi ogni giorno provare a metterle in moto tutte e farci un giro. Un bel po' di tempo, c'è il rischio che non ci resti il tempo per lavorare e pagarle queste macchine, oltre che sprecare tempo che potrebbe essere dedicato a cose più gradevoli.

Nell'ICT, e quindi nelle applicazioni crittografiche, i ragionamenti da fare sono analoghi: valutare i rischi, la probabilità che certi eventi avvengano e stimare il danno che ne deriva, quindi adottare misure ‘ragionevoli’: le misure non devono costare più del danno. Non si installa un supercomputer ridondato per essere sicuri di finire il solitario.

Va infine ricordato che la sicurezza complessiva di un sistema è data dal suo elemento più insicuro: come in una catena, per cui la resistenza è quella che deriva dall'anello più debole. E’ questo principio che sta alla base dei ragionamenti fatti al §6. Contano le chiavi, ma sono altrettanto importanti le misure con cui queste chiavi vengono protette.[66]

## Appendice 3 - Definizioni e acronimi

| | |
|---|---|
| Cleartext | Vedi Plaintext |
| attaccante | Termine generico per chiunque tenti di “rompere” un algoritmo crittografico |
| Algoritmo | un procedimento che consente di ottenere un risultato atteso eseguendo, in un determinato ordine, un insieme di passi semplici corrispondenti ad azioni scelte solitamente da un insieme finito. |
| AES | Advanced Encryption Standard |
| CESG | Communications-Electronics Security Group (UK) |
| Chiave | elemento la cui conoscenza è indispensabile per decrittare |
| Ciphertext | il testo crittografato |
| Crittoanalisi | l’arte (la scienza) di rompere un codice crittografico |
| Crittografare | trasformare un plaintext in ciphertext |
| Crittologia | l’arte (la scienza) di tenere segreto un messaggio |
| Decrittare | l’operazione inversa del crittografare |
| DES | Data Encryption Standard |
| DH | Diffie, Hellman |
| Digest | stringa di lunghezza fissa che rappresenta l’impronta di un documento |
| MCD | Massimo comune Divisore |

[66] Scegliere una password complicata per poi tenerla scritta su un *post-it* attaccato al video è un esempio classico di cattiva gestione della sicurezza.

| | |
|---|---|
| mcm | minimo comune multiplo |
| NIST | National Institute for Standard and Technology (USA) |
| NSA | National Security Agency (USA) |
| PKC | Public Key Cryptography |
| PKI | Public Key Infrastructure |
| Plaintext | il testo in chiaro |
| RFC | Request for Comment |
| RSA | Rivest, Shamir, Adlemann |
| SHA | Secure Hash Algorithm |
| Stringa | Sequenza di bit o byte |

## Bibliografia